\documentstyle[12pt]{article}
\begin{document}
\centerline{\bf Addendum to ``Defects for ample divisors of abelian
varieties,}
\centerline{\bf Schwarz lemma, and hyperbolic hypersurfaces of low degrees,''}
\centerline{\bf Amerian Journal of Mathematics 119 (1997), 1139-1172}

\bigskip
\centerline{Yum-Tong Siu and Sai-Kee Yeung\ %
\footnote{Both authors partially supported by grants from the
National Science Foundation.}
}

\bigskip
{\it This addendum was submitted to the American Journal
of Mathematics on June 6, 2000.  On the same day a copy 
was sent to J.~Noguchi who first drew our attention to a gap
in our original paper.  After three months of queries 
by e-mail, J.~Noguchi informed the first author by e-mail on 
August 28, 2000 that he finally understood the arguments.
This addendum is made available on the web at this time, 
because it has not yet appeared in a journal and is quoted in 
some recent papers of the first author.  Moreover, from
some public lectures a number of people had the wrong 
impression of the nature of the gap and did not know that it 
means only an easy modification of our original paper.}

\bigskip
J. Noguchi, J. Winkelmann, and K. Yamanoi [NWY99]
drew our attention to a difficulty in the proof of
Lemma 2 of the paper in the title,
which might arise when the projection of the space $J_k(D)$
of $k$-jets of the ample
divisor $D$ of the $n$-dimensional abelian variety $A$
completely contains $W_k$, where $W_k\subset{\bf C}^{nk}$
is the algebraic Zariski closure of the image of the differential
of the holomorphic map $f:{\bf C}\rightarrow A$.  In this addendum
we give an easy modification of our argument by
using the uniqueness part of
the fundamental theorem of ordinary differential equations.

\medskip
The modification at the same time improves somewhat
our result (see Theorem A).
For a holomorphic map to an $n$-dimensional
abelian variety $A$ whose
image is not contained in a translate of an ample divisor $D$,
the proximity function to $D$
plus the error obtained by truncating multiplicity
in the counting function at order $k_n=k_n(D)$
is of logarithmic order of
the characteristic function of the map, where $k_n$ is an explicit
function of $n$ and the Chern number $D^n$.

\medskip
The difficulty of [SY97] resulted from an attempt to use
semi-continuity of cohomology groups in deformations to avoid
employing Bloch's technique from [B26] which involves the uniqueness
part of the fundamental theorem of ordinary
differential equations.  Here we simply go back to using Bloch's
technique to show that the projection of $J_k(D)$ cannot
completely contain $W_k$ when $k$ is at least $k_n(D)$.

\medskip
A sketch of the main idea in the modification given here
is as follows.  If the local function which defines
$J_{k+1}(D)$ from $J_k(D)$ does not reduce the dimension or
the multiplicity of the part of $J_k(D)$ over $W_k$ at
a generic point, there it must
locally belong to the ideal generated by local functions defining
$J_k(D)$ and $W_k$.  The multiplicity of the part of $J_k(D)$ over
$W_k$ at a generic point
admits an effective bound depending on $D^n$ and $k$.
If the projection of $J_k(D)$ completely contains $W_k$,
then the pullback of the defining
function for a translate of $D$ by $f$ satisfies an ordinary
differential equation in some open neighborhood of $0$ in $\bf C$
whose order admits an effective bound depending on $D^n$ and $k$
and whose coefficients for the highest-order derivative is $1$.
The uniqueness part of the fundamental theorem of ordinary
differential equations yields the contradiction that some translate
of $D$ must contain the image of $f$.

\medskip
\noindent
{\it Notations and Terminology.}
Let $A$ be an abelian variety of complex dimension $n$.
We denote by ${\cal O}_A$ the structure sheaf of $A$. We denote
by $J_k(A)=A\times{\bf C}^{nk}$ the space of all $k$-jets of $A$.
We consider only subvarieties of $J_k(A)=A\times{\bf C}^{nk}$
which are algebraic along the factor ${\bf C}^{nk}$.
In other words, they are defined by local functions
on $J_k(A)=A\times{\bf C}^{nk}$ which are holomorphic
in the local coordinates of $A$ and are polynomials in
the $nk$ global coordinates of ${\bf C}^{nk}$.  Here subvarieties
of $J_k(A)=A\times{\bf C}^{nk}$ are assumed
without any further explicit mention to be
algebraic along the factor ${\bf C}^{nk}$.
Zariski closures are defined using such subvarieties.

\medskip
The projection $\pi_k:J_k(A)\rightarrow{\bf C}^{nk}$ denotes the
natural projection $J_k(A)=A\times{\bf C}^{nk}\rightarrow {\bf
C}^{nk}$ onto the second factor.  For $\ell>k$ let
$p_{k,\ell}:{\bf C}^{n\ell}\rightarrow{\bf C}^{nk}$ denote the
natural projection so that $\pi_k=p_{k,\ell}\circ\pi_\ell$.

\medskip
Let $D$ be an ample divisor in $A$.
We denote by $\theta_D$ the theta function on ${\bf C}^n$
whose divisor is $D$.  The space of $k$-jets of $D$ is denoted
by $J_k(D)$.

\medskip
For a holomorphic map $f$ from ${\bf C}$ to $A$, let $W_{k,f}$ be
the (algebraic) Zariski closure of
$\pi_k\left(\hbox{Im\,}(d^kf)\right)$ in ${\bf C}^{nk}$, where
$d^kf:J_k\left({\bf C}\right)\rightarrow J_k(A)$ is the map
induced by $f$. Clearly, $p_{k,\ell}(W_{\ell,f})=W_{k,f}$ for
$\ell>k$. When $Z_k$ is a proper subvariety of $W_{k,f}$,
$p_{k,\ell}^{-1}(Z_k)$ is a proper subvariety of $W_{\ell,f}$ for
$\ell>k$. Let $m(r,f,D)$ denote the proximity function for the map
$f$ and the divisor $D$ at the radius $r$.  Let $N_k(r,f,D)$ be
the counting function for the divisor $D$ which truncates
multiplicity at $k$, {\it i.e.,} with multiplicity replaced by $k$
if it is greater than $k$.  Let $N(r,f,D)$ denote the usual
counting function for the divisor $D$ without any truncation of
multiplicity.

\medskip
\noindent
{\it Theorem A.} Let $A$ be an abelian variety of complex dimension
$n$ and $D$ be an ample divisor
of $A$.  Inductively let $k_0=0$ and $k_1=1$ and
$$
k_{\ell+1}=k_\ell+3^{n-\ell-1}\left(4(k+1)\right)^\ell D^n
$$
for $1\leq\ell<n$.
Then for any holomorphic map $f:{\bf C}\rightarrow A$ whose image is
not contained in any translate of $D$,
$$
\displaylines{
m(r,f,D)+\left(N(r,f,D)-N_{k_n}(r,f,D)\right)
=O(\log T(r,f,D)+\log r)\qquad ||
}
$$
when $||$ means that the equation holds for $r$ outside
some set whose measure with respect to ${dr\over r}$ is
finite.

\bigskip
The rest of this note is devoted to the proof of Theorem A.
Our first step is to verify the following claim.

\bigskip
\noindent
{\it Claim 1.}  $\pi_{k_n}(J_{k_n}(D))$ does not contain $W_{k_n,f}$.

\bigskip
Assume that Claim 1 is false and we are going to derive a contradiction.

\medskip
In order to have an effective bound on the multiplicities of
the highest-dimensional branches of
$$
\pi_k^{-1}(\gamma)\cap J_k(D)=\pi_k^{-1}(\gamma)\cap\{\theta_D=d\theta_D
=\cdots=d^k\theta_D\},
$$
we introduce a mermorphic connection for the line bundle $D$ so
that the covariant differentials of $\theta_D$ with respect to
the meromorphic connection are line-bundle-valued jet differentials
on $A$, making it possible to bound the multiplicities by Chern
numbers.
By Lefschetz's theorem, $mD$ is a very ample line
bundle over $A$ for $m\geq 3$.  Choose $\tau_j\in\Gamma(A,jD)$ for
$j=3,4$ with the property that for $1\leq k\leq k_n$ there
exists a proper subvariety $Z_k$ of
$W_{k,f}$ such that
for $\gamma_k\in W_{k,f}-Z_k$ the set
$\{\tau_3\tau_4=0\}\cap J_k(D)\cap\pi_k^{-1}(\gamma_k)$
is a nowhere dense subvariety of $J_k(D)\cap\pi_k^{-1}(\gamma_k)$.

\medskip
We now define covariant differentials of $\theta_D$
by using usual differentials of meromorphic functions after
meromorphically trivialize the line bundle
$D$ by means of $\tau_3,\tau_4$.  The above condition in
the choice of $\tau_3,\tau_4$ is to make sure that the
additional pole-sets and zero-sets introduced because of
the meromorphicity of the connection are irrelevant to our
purpose.  Let
$$
{\cal D}^k\theta_D=\left(\tau_4\right)^{k+1}d^k\left(
{\theta_D\tau_3\over\tau_4}\right)
$$
for $1\leq k\leq k_n$.
Then ${\cal D}^k\theta_D$ is a $\left(4(k+1)D\right)$-valued
holomorphic $k$-jet differential on $A$ of weight $k$.
For $\ell\geq k$ and
$\gamma\in{\bf C}^{n\ell}$
we denote by $\left<{\cal D}^k\theta_D,\gamma\right>$
the element of $\Gamma(A,4(k+1)D)$ which is obtained by
evaluating ${\cal D}^k\theta_D$ at the element of
$J_k(D)$ whose image under $\pi_k$ is equal to
$p_{k,\ell}(\gamma)$.

\medskip
We use the convention that
a set of complex dimension $-1$ means the empty set.

\medskip
For $1\leq\ell\leq n$ we introduce the following statement
$(A)_\ell$ and use the convention that $(A)_0$ always holds.
Let $0\leq\ell_0\leq n$ be the largest integer
inclusively between $0$ and $n$
such that the following statement $(A)_\ell$ holds
with $\ell=\ell_0$.

\bigskip
\noindent
$(A)_\ell$.  There exists a proper subvariety $E_\ell$ of $W_{k_\ell,f}$
such that the complex dimension of the common zero-set of
$$
\theta_D,\left<{\cal D}\theta_D,\gamma\right>,
\left<{\cal D}^2\theta_D,\gamma\right>,\cdots,
\left<{\cal D}^{k_\ell}\theta_D,\gamma\right>
$$
in $A-\left\{\tau_3\tau_4=0\right\}$ is no more than $n-1-\ell$
for $\gamma\in W_{k_\ell,f}-E_\ell$.

\bigskip
We can choose

\medskip
\noindent
(i) $\sigma_0\in\Gamma\left(A,4(k_{\ell_0}+1)D\right)$
and $\sigma_j\in\Gamma\left(A,4(k_{\ell_0}-j)D\right)$
for $1\leq j<k_{\ell_0}$,

\medskip
\noindent
(ii) complex numbers $c_{i,j}$ for $1\leq i\leq\ell_0$ and
$1\leq j\leq k_{\ell_0}$, and

\medskip
\noindent (iii) a proper subvariety $\tilde E_{\ell_0}$ of
$W_{k_{\ell_0},f}$ containing $E_{\ell_0}$,

\medskip
\noindent
such that, when we set $\sigma_{k_{\ell_0}}\equiv 1$ and
$$
\rho_{i,\gamma}=\sum_{j=0}^{k_{\ell_0}}c_{i,j}
\sigma_j\left<{\cal D}^j\theta_D,\gamma\right>
$$
for $1\leq i\leq\ell_0$ and $\gamma\in W_{k_{\ell_0},f}
-\tilde E_{\ell_0}$, the common zero-set of
$$
\rho_{1,\gamma},\cdots,\rho_{\ell_0,\gamma}
\in\Gamma\left(A,4(k_{\ell_0}+1)D\right)
$$
in $D-\left\{\tau_3\tau_4=0\right\}$ is a subvariety of pure
complex dimension $n-1-\ell_0$ in $D$.

\bigskip Because of the assumption that
$\pi_{k_n}\left(J_{k_n}(D)\right)$ contains $W_{k_n,f}$, we know
that $\ell_0<n$. Let $\kappa_n=\theta_D$.  We can choose
$$
\kappa_1,\cdots,\kappa_{n-{\ell_0}-1}\in\Gamma(A,3D)
$$
such that the intersection of the common zero-set of
$\kappa_1,\cdots,\kappa_{n-\ell_0-1}$ and every
$\ell_0$-dimensional branch of the common zero-set of
$\rho_{1,\gamma},\cdots,\rho_{\ell_0,\gamma}$ in
$D-\left\{\tau_3\tau_4=0\right\}$ is zero-dimensional. Then the
sum, over all points $Q\in A-\left\{\tau_3\tau_4=0\right\}$, of
the dimension over $\bf C$ of the (vector-space) stalk at $Q$ of
the sheaf
$$
{\cal O}_A\,\bigg/
\left(\sum_{j=1}^{n-\ell_0}{\cal O}_A\kappa_j+
\sum_{i=1}^{\ell_0}{\cal O}_A\rho_{i,\gamma}\right)
$$
is no more than
$$
3^{n-\ell_0-1}\left(4(k_{\ell_0}+1)\right)^{\ell_0}D^n.
$$
Since the above estimates imply that the complex dimension $d_s$
of the stalk at $Q\in D-\left\{\tau_3\tau_4=0\right\}$ of
$$
{\cal O}_A\,\bigg/ \left(\sum_{j=1}^{n-\ell_0}{\cal O}_A\kappa_j+
\sum_{i=1}^s{\cal O}_A \left<{\cal
D}^s\theta_D,\gamma\right>\right)
$$
is no more than $k_{\ell_0+1}-k_{\ell_0}$ for
$\gamma\in\left(p_{k_{\ell_0},k_{\ell_0+1}}\right)^{-1}
\left(W_{k_{\ell_0},f}-\tilde E_{\ell_0}\right)$ and
$k_{\ell_0}\leq s\leq k_{\ell_0+1}$, it follows that there exists
some $k_{\ell_0}< q\leq k_{\ell_0+1}$ such that $d_{q-1}= d_q$.
The $q$ depends on the choice of $
Q,\gamma,\kappa_1,\cdots,\kappa_{n-\ell_0}$, but we can make $q$
independent of the choice for $Q$ in some Zariski open subset of
$D-\left\{\tau_3\tau_4=0\right\}$, for
$\left(\kappa_1,\cdots,\kappa_{n-\ell_0-1}\right)$ in some Zariski
open subset of the product of $n-\ell_0-1$ copies of
$\Gamma\left(A,3D\right)$, and for $\gamma$ in some Zariski open
subset of $\left(p_{k_{\ell_0},k_{\ell_0+1}}\right)^{-1}
\left(W_{k_{\ell_0},f}-\tilde E_{\ell_0}\right)$.

\bigskip Since $(A)_{\ell_0+1}$ does not hold, there exists a proper
subvariety $F_{\ell_0+1}$ of $W_{\ell_0+1,f}$ containing
$p_{\ell_0,\ell_0+1}^{-1}\left(\tilde E_{\ell_0}\right)$ such that
that the complex dimension of the common zero-set of
$$
\theta_D,\left<{\cal D}\theta_D,\gamma\right>, \left<{\cal
D}^2\theta_D,\gamma\right>,\cdots, \left<{\cal
D}^{k_{\ell_0+1}}\theta_D,\gamma\right>
$$
in $A-\left\{\tau_3\tau_4=0\right\}$ is $n-1-\ell_0$ for
$\gamma\in W_{k_{\ell_0+1},f}-F_{\ell_0+1}$.  Let $V_{\ell_0+1}$
be the minimum subvariety in
$J_q(D)\cap\pi_q^{-1}\left(W_{q,f}\right)$ which contains all the
branches of dimension $<n-1-\ell_0$ of
$J_q(D)\cap\pi_q^{-1}(\gamma)$ for all $\gamma\in W_{q,f}$.  Then
$V_{\ell_0+1}$ is a proper subvariety of
$J_q(D)\cap\pi_q^{-1}\left(W_{q,f}\right)$.  After replacing
$V_{\ell_0+1}$ by a larger proper subvariety of
$J_q(D)\cap\pi_q^{-1}\left(W_{q,f}\right)$ to take care of the
requirement of $Q,\gamma,\kappa_1,\cdots,\kappa_{n-\ell_0-1}$
being in some appropriate Zariski open subsets in the choice of
$q$ made above, we conclude that the following statement
$(B)_{\ell}$ holds with $\ell=\ell_0+1$.

\bigskip
\noindent $(B)_\ell$. There exist $k_{\ell-1}<q\leq k_\ell$, a
subvariety $V_\ell$ of $J_q(D)\cap\pi_q^{-1}(W_{q,f})$ containing
$J_q(D)\cap\left\{\tau_3\tau_4=0\right\}$, and a proper subvariety
$F_\ell$ of $W_{q,f}$, such that

\medskip
\noindent
(i) $\left(J_q(D)-V_\ell\right)\cap\pi_q^{-1}(\gamma_q)$ is a
nonempty $(n-\ell)$-dimensional subvariety of $J_q(A)-V_\ell$ for every
$\gamma_q\in W_{q,f}-F_\ell$, and

\medskip
\noindent
(ii) on $J_{q}(A)-\pi_q^{-1}\left(F_\ell\right)-V_\ell$
the function ${\cal D}^q\theta_D$ locally belongs to the ideal
sheaf which is
locally generated by the functions
$$
\theta_D,{\cal D}\theta_D,{\cal D}^2\theta_D,
\cdots,{\cal D}^{q-1}\theta_D
$$
and by the pullbacks to $J_q(A)$ of the local holomorphic
functions on ${\bf C}^{nq}$ vanishing on $W_{q,f}$.

\bigskip
The numbers $k_0=0$ and $k_1=1$ are chosen so that $(B)_\ell$ holds
with $\ell=1$ when $(A)_\ell$ fails for $\ell=1$.

\medskip
There exists some $\zeta_0\in{\bf C}$ such that
$$
\pi_q\left((d^q f)(\zeta_0)\right)\in W_{q,f}-F_{\ell_0+1},
$$
where $\left(d^qf\right)(\zeta_0)$ denotes the element of $J_q(A)$
defined by $f$ at $\zeta_0$. By Condition (i) of $(B)_{\ell_0+1}$
there exists some point $Q_0\in D$ such that when we define
$$
g(\zeta)=\left(Q_0-f(\zeta_0)\right)+f(\zeta)\in A
$$
for $\zeta\in{\bf C}$, we have
$$
\left(d^q g\right)\left(\zeta_0\right)\in
J_q(D)-V_{\ell_0+1}-\pi_q^{-1}\left(F_{\ell_0+1}\right).\leqno{(1)}
$$
Since $\pi_q\left(\left(d^q
g\right)\left(\zeta\right)\right)=\pi_q\left(\left(d^q
f\right)\left(\zeta\right)\right)\in W_{q,f}$ for $\zeta\in{\bf
C}$, from Condition (ii) of $(B)_{\ell_0+1}$ it follows that
$\theta_D\circ g$ satisfies a differential equation
$$
\left({d^q\over d\zeta^q}\right)(\theta_D\circ g)(\zeta)
=\sum_{j=0}^{q-1}h_j(\zeta)
\left({d^j\over d\zeta^j}\right)(\theta_D\circ g)(\zeta)
\leqno{(2)}
$$
on some open neighborhood $U$ of $\zeta_0$ in ${\bf C}$ for
some holomorphic functions $h_0(\zeta),\cdots,h_{q-1}(\zeta)$
on $U$.  By $(1)$ and $(2)$
and the uniqueness part of the
fundamental theorem of ordinary differential equations,
$(\theta_D\circ g)(\zeta)$ is identically zero for $\zeta\in{\bf C}$.
Hence the image of $f$ is contained in the translate of $D$
by $Q_0-f(\zeta_0)$. This contradicts the assumption that
the image of $f$ is not contained in any translate of $D$
and therefore finishes the verification of Claim 1.

\bigskip
Claim 1 yields the following lemma.

\medskip
\noindent
{\it Lemma 1.}
There exists a polynomial
$$
P=P\left(\{d^\mu z_\nu\}_{1\leq\mu\leq k_n,1\leq\nu\leq n}\right)
$$
in the variables
$d^\mu z_\nu$ ($1\leq\mu\leq k_n,1\leq\nu\leq n$) with
coefficients in $\bf C$ such that
$P$ is not identically zero on $W_{k_n,f}$ and
$\pi_{k_n}^*P$ is identically zero on
$J_{k_n}(D)$.

\medskip
Let $h_{\alpha,\bar\beta}$ be the positive definite Hermitian matrix
such that
$$
\sum_{\alpha,\beta=1}^n
h_{\alpha,\bar\beta}dz_\alpha\overline{dz_\beta}
$$
represents the Chern class of the ample line bundle $D$.
Denote by $\nabla$
the covariant differentiation for sections of $D$ with respect to the metric
of $D$ whose curvature form is
$\sum_{\alpha,\beta=1}^n
h_{\alpha,\bar\beta}dz_\alpha\overline{dz_\beta}$.
Then for a local section $s$ of the ample line bundle $D$,
$$\nabla s
=ds+\left(\sum_{\alpha,\beta=1}^n
h_{\alpha,\bar\beta}\overline{z_\beta}
dz_\alpha\right)s.
$$

\medskip
The following lemma is an immediate consequence of Lemma 1.

\medskip
\noindent
{\it Lemma 2}.
The pullback $\tilde P$ to ${\bf C}^n$ of $P$ can be written in the form
$$
\tilde P=\sum_{\mu=0}^{k_n}\tilde\rho_\mu
\left(\nabla^\mu\theta_D\right),\leqno{(3)}
$$
where
$\tilde\rho_\mu$ is a polynomial in $d^\mu z_\nu$
($1\leq\mu\leq k_n,1\leq\nu\leq n$) whose coefficients are
bounded smooth functions on ${\bf C}^n$.

\medskip
Let ${\cal A}_r(\cdot)$ denote the average over the circle of
radius $r$ centered at $0$.
Let $T(r,f,D)$ be the Nevanlinna characteristic function for the
function $f$ and the ample line bundle $D$ of $A$ at the radius $r$.
Let $\tilde f:{\bf C}\rightarrow{\bf C}^n$ be the lifting of
$f:{\bf C}\rightarrow A$.  For a meromorphic function $F$ on ${\bf C}$
we denote by $N(r,F,\infty)$ the counting function of the map
${\bf C}\rightarrow{\bf P}_1$ defined by $F$ and denote by
$T(r,F)$ its Nevanlinna characteristic function (with respect to
the line bundle of ${\bf P}_1$ of degree $1$).

\medskip
\noindent
{\it Final Step in the Proof of Theorem A.}
From $d^\mu z_\nu=d^\mu\log\left(e^{z_\nu}\right)$ and
the logarithmic derivative lemma it follows that
$$
{\cal A}_r\left(\log^+\left|\tilde f^*
\left(d^\mu z_\nu\right)\right|\right)
=O(\log T(r,f,D)+\log r)\quad ||
$$
and
$$
{\cal A}_r\left(\log^+\left|\tilde f^*
\tilde P\right|\right)
=O(\log T(r,f,D)+\log r)\quad ||.\leqno{(4)}
$$
Here we interpret $\tilde f^*d^\mu z_\nu$ as
a function of $\zeta$ which is the coefficient of the power of
$d\zeta$.  The notation $\tilde f^*\tilde P$ above
and $\tilde f^*\left({\tilde P\over\theta_D}\right)$
carry similar meanings and are functions of $\zeta$.
By (3) and the logarithmic derivative lemma it follows that
$$
{\cal A}_r\left(\log^+\left|\tilde f^*
\left({\tilde P\over\theta_D}\right)\right|\right)
=O(\log T(r,f,D)+\log r)\quad ||.\leqno{(5)}
$$
From (3) we conclude that
$$
N\left(r,\tilde f^*\left({\tilde P\over\theta_D}\right),\infty\right)
\leq N_{k_n}(r,f,D).\leqno{(6)}
$$
It follows from the First Main Theorem of Nevanlinna theory that
$$
\displaylines{
m(r,f,D)+N(r,f,D)\cr
=T\left(r,\tilde f^*\left({1\over\theta_D}\right)\right)
+O(1)\cr
=T\left(r,\tilde f^*\left({\tilde P\over\theta_D}
{1\over\tilde P}\right)\right)+O(1)\cr
\leq T\left(r,\tilde f^*\left({\tilde P\over\theta_D}
\right)\right)+
T\left(r,{1\over\tilde P}\right)+O(1).\cr
}
$$
Using (4), (5), and (6), we obtain
$$
m(r,f,D)+N(r,f,D)
\leq
N_{k_n}(r,f,D)+O\left(\log T\left(r,f,D\right)+\log r\right)\quad ||.
$$
Q.E.D.

\medskip
\noindent
{\it Remark.} The same arguments work for semi-abelian varieties
with straightforward modifications.

\bigskip
\noindent
{\it References.}

\medskip
\noindent
[B26] A. Bloch, Sur les syst\`emes de fonctions uniformes satisfaisant
\`a l'\'equa\-tion d'une vari\'t\'e alg\'ebrique dont l'irr\'egularit\'e
d\'epasse la dimension, J. de Math. 5 (1926), 19-66.

\medskip
\noindent
[NWY99] J. Noguchi, J. Winkelmann, K. Yamanoi, The second main
theorem for holomorphic curves into semi-Abelian varieties, Preprint,
University of Tokyo, 1999.

\medskip
\noindent
[SY97] Y.-T. Siu and S.-K. Yeung,
Defects for ample divisors of abelian varieties, Schwarz lemma, and
hyperbolic hypersurfaces of low degrees, Amer. J.
Math. 119 (1997), 1139-1172.

\bigskip
\noindent
{\it Authors' addresses:}

\medskip
\noindent
Yum-Tong Siu, Department of Mathematics, Harvard University, Cambridge,
MA 02138. {\it e-mail:} siu@math.harvard.edu.

\medskip
\noindent
Sai-Kee Yeung, Department of Mathematics, Purdue University, West Lafay\-ette,
IN 47907. {\it e-mail:} yeung@math.purdue.edu.

\end{document}